\documentclass[journal]{IEEEtran}
\usepackage{epsfig,graphicx,amssymb,amsmath}
\usepackage{bm}
\usepackage{amsthm}
\usepackage{amsmath,bm}
\usepackage{subfigure}
\usepackage{algorithm, algorithmic}
\usepackage{amsfonts}
\usepackage[noadjust]{cite}  %% for citation compress
\usepackage{url}
\usepackage{array,tabularx}
\usepackage{indentfirst}
\usepackage{color}
\usepackage{stfloats}

\ifCLASSINFOpdf
\else
\fi
\begin{document}
\title{{A Unified Distributed Method for Constrained Networked Optimization via Saddle-Point Dynamics}}
%


\author{Yi Huang, Ziyang Meng, \IEEEmembership{Senior Member, IEEE}, Jian Sun, \IEEEmembership{Senior Member, IEEE}, and Wei Ren,  \IEEEmembership{Fellow, IEEE}% <-this % stops a space
\thanks{This work has been supported in part by the National
Natural Science Foundation of China under Grants under Grants 62103223, 61925303, 62088101, 61833009 and U19B2029.

Yi Huang and Jian Sun are with the School of Automation, Beijing Institute of Technology, Beijing 100081, China (e-mail: yihuang@bit.edu.cn, e-mail: sunjian@bit.edu.cn).

Ziyang Meng is with the Department of Precision
Instrument, Tsinghua University, Beijing 100084, China (e-mail: ziyangmeng@mail.tsinghua.edu.cn).

Wei Ren is with Department of Electrical and Computer Engineering, University of California, Riverside, CA
92521, USA (e-mail: ren@ee.ucr.edu).
}}% <-this % stops a space
\maketitle

\begin{abstract}
This paper develops a unified distributed method for solving two classes of constrained networked optimization problems, i.e., optimal consensus problem and resource allocation problem with non-identical set constraints. We first transform these two constrained networked optimization problems into a unified saddle-point problem framework with set constraints. Subsequently, two projection-based primal-dual algorithms via Optimistic Gradient Descent Ascent (OGDA) method and Extra-gradient (EG) method are developed for solving constrained saddle-point problems. It is shown that the developed algorithms achieve exact convergence to a saddle point with an ergodic convergence rate $O(1/k)$ for general convex-concave functions. Based on the proposed primal-dual algorithms via saddle-point dynamics, we develop unified distributed algorithm design and convergence analysis for these two networked optimization problems. Finally, two numerical examples are presented to demonstrate the theoretical results.
\end{abstract}


\newtheorem{Assumption}{Assumption}
 \newtheorem{Remark}{Remark}
 \newtheorem{Lemma}{Lemma}
 \newtheorem{Definition}{Definition}
  \newtheorem{Proposition}{Proposition}
 \newtheorem{Theorem}{Theorem}
 \newtheorem{Property}{Property}
 \newtheorem{Corollary}{Corollary}
 \newtheorem{Example}{Example}

\section{Introduction}

The problem of distributed optimization has attracted considerable attention in recent decades due to its wide applications in machine learning, power systems, multi-robot localization, sensor networks, and resource allocation \cite{3}. In general, most distributed optimization problems in the existing literature can be divided into two categories: optimal consensus problem and optimal resource allocation problem \cite{4a}. The main difference of these two problems is that in the first problem, each agent has its own objective function with respect to a common decision variable, while in the second one, all the agents own independent local objective functions and decision variables but these decision variables are coupled in a global equality constraint. To solve the optimal consensus problem, a common approach is to introduce a consensus constraint such that the coupled objective functions can be separated with the local decision variables. In such a case, the optimal consensus problem and optimal resource allocation problem can be both regarded as a class of optimization problems with a linear equality constraint. For these two classes of optimization problems, many discrete-time and continuous-time algorithms are developed in \cite{4a, 5,6,7,7a,8,9,10,11,12,14,14a,14b,15,16,17,17a,17b,17c}.
%
%Distributed optimization problem has attracted considerable attention due to its wide applications in machine learning, power system, multi-robot localization, sensor network, and resource allocation \cite{1,2,3,4}. Most distributed optimization problems in the existing literatures can be divided into two categories. The first one is the optimal consensus problem, i.e., \emph{the decision-coupled problem}, in which the local objective functions of all the agents are coupled with respect to a common decision variable. The second one is the optimal resource allocation problem, i.e., \emph{the constrained-coupled problem}, in which each agent has local objective function while all the agents' local decision variables are coupled with respect to a linear equality constraint. To solve the optimal consensus problem, a common approach is to introduce a consensus constraint such that the coupled objective functions can be separated with the local decision variables. In such a case, the optimal consensus problem and optimal resource allocation problem can be both regarded as a class optimization problem with a linear equality constraint. For these two optimization problems, many discrete-time and continuous-time algorithms are developed in \cite{5,6,7,8,9,10,11,12,13,14,15,16,17}.

Note that most existing distributed optimization algorithms to solve the optimal consensus problem and resource allocation problem are designed separately. Fewer results provide a unified framework for analysis and design of these two optimization problems. As mentioned above, these two optimization problems can be both viewed as the constrained optimization problem with a linear equality constraint. For the constrained optimization problem, we can transform it to a class of saddle-point problems in terms of the corresponding Lagrangian functions \cite{18a}. This fact illustrates that the above-mentioned two optimization problems can both be transformed into the saddle-point problems. Therefore, when the saddle points of the corresponding Lagrangian functions are obtained, these two optimization problems can be solved.

It is well known that saddle-point problems arise in many areas such as constrained optimization \cite{18b}, robust control \cite{18c}, zero-sum games \cite{18} and generative adversarial networks (GANs) \cite{19a}. Some typical first-order optimization methods (e.g., Gradient Descent Ascent (GDA), Optimistic Gradient Descent Ascent (OGDA) and Extra-gradient (EG) methods) have been proposed to solve the saddle-point problems. This paper focuses on OGDA and EG methods, whose ideas were first proposed in \cite{20} and \cite{21}, and have attracted considerable attention. The authors of \cite{23} showed the linear convergence rates of OGDA and EG methods for a special case, i.e., $f(x,y)=x^{T}Ay$, where $A$ is square and full rank. In \cite{24}, the authors proposed a variant of EG method with linear convergence when $f(x,y)$ is strongly convex-strongly concave, and applied it to the GANs training. The authors of \cite{25} showed OGDA and EG methods as approximate variants of the proximal point method, and provide their linear convergence for strongly convex-strongly concave functions. For general convex-concave functions, the authors of \cite{26} provided a unified convergence analysis of OGDA and EG methods and proved that these two methods can both achieve an ergodic convergence rate of $O(1/k)$. Nevertheless, the last iteration of \cite{26} is shown to only converge into a bounded neighborhood of a saddle point instead of achieving exact convergence to a saddle point. In addition, we note that most results on OGDA and EG methods mentioned above only consider the saddle-point problems in absence of constraints. Actually, the saddle-point problems with set constraints are very common in practical applications.

Inspired by the above discussions, this paper tries to establish the relationship between two classes of constrained networked optimization problems and general constrained saddle-point problems, and then solve them under a unified saddle-point dynamics framework. Compared with the related results, the main contributions of this paper are three-fold.
\begin{enumerate}
  \item [\textbf{c1)}] We develop unified distributed algorithm design and convergence analysis via saddle-point dynamics to solve two classes of constrained networked optimization problems, i.e., optimal consensus problem and resource allocation problem with non-identical set constraints.

  \item [\textbf{c2)}] Two projection-based primal-dual algorithms via OGDA and EG methods are developed for constrained saddle-point problem for general convex-convex functions. Unlike the results of \cite{26} that are only shown to converge into a bounded neighborhood of a saddle point, the developed algorithms achieve exact convergence to a saddle point with an ergodic convergence rate $O(1/k)$.

  \item [\textbf{c3)}] The developed distributed algorithms are with constant step-sizes and performs better convergence performance than the algorithms in \cite{14,15,16} with diminishing step-sizes. In contrast with the constant step-size algorithms in \cite{4a} and \cite{14a}, the developed algorithms are more easier to be implemented without solving the sub-optimization problem at each iteration.

\end{enumerate}

%Also, compared with gradient tracking algorithm of \cite{14b}, the developed distributed algorithms effectively solve the optimal consensus problem with non-identical set constraints.

The rest of this paper is organized as follows. Section II formulates the considered problem. Section III proposes two primal-dual algorithms via OGDA and EG methods. Section IV develops unified distributed algorithms to solve two networked optimization problems. Section V gives simulation examples and Section VI concludes this paper.

\section{Preliminaries and formulation}
Notation: Let $\mathbb{R}$ be the set of real numbers and $\mathbb{N}$ be the set of natural numbers. $I_{n}$ is the $n\times n$ identity matrix and $\bm 1_{n}$ is the $n\times 1$ ones vector. $\Vert \cdot\Vert$ denotes the Euclidean norm. Let $\mathcal{I}_{N}=\{1,2,\ldots, N\}$ and  $\text{col}(x_{i})^{N}_{i=1}$ be a column stack of the vector $x_{i}, i\in \mathcal{I}_{N}$. $\text{diag}(W_{i})_{i=1}^N$ denotes a diagonal block matrix and $W_{i}$ is placed in the $i$th diagonal block, and $\otimes$ represents the Kronecker product.

\subsection{Problem Formulation}
In this section, two classes of constrained networked optimization problems are formulated. Consider a network graph $\mathcal{G}$ of $N$ agents. The distributed optimal consensus problem with non-identical set constraints is described by \cite{5}
\begin{align}\label{2}
 \min \limits_{x \in \Omega} f(x)=\sum^{N}_{i=1} f_{i}(x_{i}), ~\text{s.t.}~(L\otimes I_{m})x=0,
\end{align}
where $x_{i}\in \mathbb{R}^{m}, x=\text{col}(x_{i})^{N}_{i=1} \in \mathbb{R}^{Nm}$, $\Omega=\prod^{N}_{i=1}\Omega_{i}$ is the Cartesian product, and $L\in \mathbb{R}^{N\times N}$ is an Laplacian matrix of graph $\mathcal{G}$. In this problem, each agent only privately has access to local objective function $f_{i}(x_{i})$ and set constraint $\Omega_{i}, i\in \mathcal{I}_{N}$. Provided that graph $\mathcal{G}$ is connected, $(L\otimes I_{m})x=0$ implies that the consensus $x_{i}=x_{j}$ is satisfied for $\forall i,j\in \mathcal{I}_{N}$.

Next, the distributed resource allocation problem via a multi-agent network is formulated as \cite{11}
\begin{align}\label{3}
\min \limits_{y\in \Omega} h(y)=\sum^{N}_{i=1}h_{i}(y_{i}), ~\text{s.t.}~\sum^{N}_{i=1}W_{i}y_{i}=\sum^{N}_{i=1}d_{i},
\end{align}
where $h_{i}(y_{i}): \mathbb{R}^{q_{i}} \to \mathbb{R}$ is the local objective function of agent $i$, $y_{i}\in \mathbb{R}^{q_{i}}$ is its local decision variable, $y=\text{col}(y_{i})^{N}_{i=1}\in \mathbb{R}^{q}$ with $q=\sum^{N}_{i=1}q_{i}$, and $\Omega=\prod^{N}_{i=1}\Omega_{i}$ is the Cartesian product. $\sum^{N}_{i=1}W_{i}y_{i}=\sum^{N}_{i=1}d_{i}$ is the coupled equality constraint, in which $W_{i} \in \mathbb{R}^{m\times q_{i}}$ and $d_{i}\in \mathbb{R}^{m}$ are the local data only known by agent $i$.

The following standard assumptions are imposed.
\begin{Assumption}
(i) The graph $\mathcal{G}$ is undirected and connected. (ii) $\Omega_{i}, i\in \mathcal{I}_{p}$, is closed and convex. The local objective functions $f_{i}(x_{i})$ and $h_{i}(y_{i})$ are differentiable and convex on $\Omega_{i}$, and their gradients $\nabla f_{i}(x_{i})$ and $\nabla h_{i}(y_{i})$ are Lipschitz continuous for $\forall i\in \mathcal{I}_{N}$. (iii) There exists at least one solution to the problems \eqref{2} and \eqref{3}.
\end{Assumption}

\begin{Remark}
The problems \eqref{2} and \eqref{3} capture a wide class of networked optimization problems in practical applications. For instance, the optimal rendezvous, cooperative localization and machine learning in \cite{3} can be described by problem \eqref{2}. The resource scheduling, economic dispatch and flow control in smart grids \cite{10,11,12} can be formulated by problem \eqref{3}.
\end{Remark}

We establish the relationships between the above two classes of constrained networked optimization problems and constrained saddle-point problems for general convex-concave functions. For the problem \eqref{2}, its augmented Lagrangian function is $L_{1}(x,v)=\sum^{N}_{i=1} f_{i}(x_{i})+v^{T}(L\otimes I_{m})x+\frac{1}{2}x^{T}(L\otimes I_{m})x$, where $v=\text{col}(v_{i})^{N}_{i=1}\in \mathbb{R}^{Nm}$ is the dual variable \cite{7}. Then, the optimization problem \eqref{2} can be transformed into the following constrained saddle-point problem
\begin{align}\label{2a}
\min_{x \in \Omega}\max_{v\in \mathbb{R}^{Nm}} L_{1}(x,v).
\end{align}
Note that $L_{1}(x,v)$ is a convex-concave function. We have that the problem \eqref{2} is reformulated as a constrained saddle-point problem for general convex-concave functions.

For the problem \eqref{3}, its modified Lagrangian function can be derived as $L_{2}(y,z,\lambda)=\sum^{N}_{i=1}h_{i}(y_{i})+\lambda^{T}(Wy-d-(L\otimes I_{m})z)-\frac{1}{2}\lambda^{T}(L\otimes I_{m})\lambda$,
where $W=\text{diag}(W_{i})^{N}_{i=1}\in \mathbb{R}^{Nm\times q}, d=\text{col}(d_{i})^{N}_{i=1}\in \mathbb{R}^{Nm}$, $\lambda=\text{col}(\lambda_{i})^{N}_{i=1}\in \mathbb{R}^{Nm}$ is the dual variable, and $z=\text{col}(z_{i})^{N}_{i=1}\in \mathbb{R}^{Nm}$ is an auxiliary variable (see eq. (12) in \cite{11}). The problem \eqref{3} is transformed into the following constrained saddle-point problem \cite{11}
\begin{align}\label{4}
\min_{y \in \Omega, z\in \mathbb{R}^{Nm}}\max_{\lambda \in \mathbb{R}^{Nm}} L_{2}(y,z,\lambda).
\end{align}
Similarly, we have that $L_{2}(y,z,\lambda)$ is a convex-concave function. This implies that problem \eqref{3} can be also transformed into a constrained saddle-point problem for general convex-concave functions.

\subsection{Unified Problem Framework}
To solve the above two classes of constrained networked optimization problems via a unified framework, we consider the following general constrained saddle-point problem
\begin{align}\label{1}
\min_{x \in \mathcal{X}}\max_{y\in \mathcal{Y}} f(x,y),
\end{align}
where $\mathcal{X} \subseteq \mathbb{R}^{n}$ and $\mathcal{Y}\subseteq \mathbb{R}^{m}$ are both closed and convex, and $f: \mathcal{X}\times \mathcal{Y}\to \mathbb{R}$ is a convex-concave objective function, i.e., for any $y\in \mathcal{Y}$, $f(x,y)$ is a convex function with respect to $x\in \mathcal{X}$, and for any $x \in \mathcal{X}$, $f(x,y)$ is a concave function with respect to $y\in \mathcal{Y}$. We focus on finding a saddle point $(x^{*},y^{*})\in \mathcal{X}\times \mathcal{Y} $ of problem \eqref{1} that satisfies $f(x^{*},y)\le f(x^{*},y^{*})\le f(x,y^{*}), \forall (x,y)\in \mathcal{X}\times \mathcal{Y}$.

According to the optimal condition of \cite{18b}, the pair $(x^{*},y^{*})$ is a saddle point of \eqref{1} if the following variational inequality holds for $\forall (x,y)\in \mathcal{X}\times \mathcal{Y}$.
\begin{align}\label{1a}
\begin{bmatrix} \nabla_{x} f(x^{*},y^{*})\\ -\nabla_{y} f(x^{*},y^{*})\end{bmatrix}^{T}\begin{bmatrix} x-x^{*} \\ y-y^{*}\end{bmatrix} \ge 0.
\end{align}
%\begin{Remark}
%Problem \eqref{1} is a widely investigated model that appears in many applications, including game theory, GAN, multi-agent reinforcement learning and statistical learning \cite{25,26}.
%\end{Remark}
\begin{Assumption}
The function $f(x,y)$ is continuously differentiable for any $x\in \mathcal{X}$ and $y\in \mathcal{Y}$. The gradient $\nabla_{x} f(x,y)$ is $l_{xx}$-Lipschitz in $x$, and $l_{xy}$-Lipschitz in $y$. The gradient $\nabla_{y} f(x,y)$ is $l_{yx}$-Lipschitz in $x$, and $l_{yy}$-Lipschitz in $y$. If $f(x,y)=x^{T}By$ is a bilinear function with constant matrix $B$, we obtain that the Lpschitz constants $l_{xx}$ and $l_{yy}$ are zero.
\end{Assumption}

\begin{Assumption}
The solution set of problem \eqref{1} is nonempty.
\end{Assumption}

This paper aims to develop a unified distributed method for solving two classes of networked optimization problems. To achieve this goal, we first propose two primal-dual algorithms for constrained saddle-point problem \eqref{1}, and then develop unified distributed algorithms via saddle-point dynamics for constrained networked optimization problems \eqref{2} and \eqref{3}.
%The aim of this paper is to design a distributed control input $u_{i}$ such that all the agents can cooperatively seek a common point $x^{*}_{0}$ that minimizes the sum of local objective functions, i.e.,
%\begin{align}
%x^{*}_{0} \in
%\end{align}

\section{Saddle-Point Dynamics Design}
In this section, we first develop two projection-based primal-dual algorithms by using OGDA and EG methods to solve the constrained saddle-point problem \eqref{1}. Next, the convergence analysis of these two algorithms is provided.

\subsection{Primal-dual algorithm via OGDA method}
We develop a projection-based primal-dual algorithm via OGDA to solve the constrained saddle-point problem \eqref{1}
\begin{subequations}\label{d1}
\begin{align}
x_{k+1}&=\mathcal{P}_{\mathcal{X}}\Big(x_{k}-\alpha \nabla_{x} f(x_{k},y_{k}) \\
&\quad-\alpha(\nabla_{x} f(x_{k},y_{k})-\nabla_{x} f(x_{k-1},y_{k-1})\Big ), \notag \\
y_{k+1}&=\mathcal{P}_{\mathcal{Y}}\Big(y_{k}+\alpha \nabla_{y} f(x_{k},y_{k})\\
&\quad+\alpha(\nabla_{y} f(x_{k},y_{k})-\nabla_{y} f(x_{k-1},y_{k-1}))\Big), \notag
\end{align}
\end{subequations}
where $\mathcal{P}_{\mathcal{X}}(\cdot)$ and $\mathcal{P}_{\mathcal{Y}}(\cdot)$ represent the projection operations on $\mathcal{X}$ and $\mathcal{Y}$, respectively, and $\alpha$ is the constant step-size that will be specified later.

Let $z=\text{col}(x,y)\in \mathcal{X}\times \mathcal{Y} \subset \mathbb{R}^{m+n}$ and the operator $F: \mathcal{X}\times \mathcal{Y} \to \mathbb{R}^{m+n}$ as $F(z)=\text{col}(\nabla_{x} f(x,y), -\nabla_{y} f(x,y))$. Eq. \eqref{d1} can be arranged as
\begin{align}\label{d2}
z_{k+1}=\mathcal{P}_{\Lambda}\Big(z_{k}-\alpha F(z_{k})-\alpha(F(z_{k})-F(z_{k-1}))\Big),
\end{align}
where $\Lambda=\mathcal{X}\times \mathcal{Y}$.

In contrast to the GDA algorithm that is formulated as $z_{k+1}=\mathcal{P}_{\Lambda}(z_{k}-\alpha F(z_{k}))$ in \cite{27}, the main difference of the proposed OGDA-based algorithm \eqref{d2} is the added gradient correction term $-\alpha (F(z_{k})-F(z_{k-1}))$, which includes the gradient information of $f(x,y)$ at the current iteration and previous iteration. The advantages of adding the gradient correction term is to guarantee exact convergence to a saddle point of general convex-concave functions. As mentioned in \cite{28}, the GDA algorithm of \cite{27} requires strongly convex-strongly concave condition of objective function to ensure the exact convergence and may not converge to a saddle point for general convex-concave functions. This result is also illustrated in Example 1 given in the following simulation section.

%\begin{Remark}
%Inspired by the work in \cite{26}, it is not hard to show that the proposed OGDA-based algorithm \eqref{d2} is regarded as an approximation of the proximal point approach (PPA) for the constrained saddle-point problem \eqref{1}, i.e.,
%\begin{align}\label{d2a}
%z_{k+1}=\mathcal{P}_{\Lambda}\Big(z_{k}-\alpha F(z_{k+1})+\varepsilon_{k}\Big),
%\end{align}
%with $\varepsilon_{k}=\alpha \{F(z_{k+1})-F(z_{k})-(F(z_{k})-F(z_{k-1}))\}$. Note that if $\varepsilon_{k}=0$, the algorithm \eqref{d2} reduces to the PPA for the problem \eqref{1}. Note that PPA achieves a sublinear convergence of $O(1/k)$ for general convex-concave functions. Then, a conjecture is that the developed OGDA-based algorithm \eqref{d2a} has the same convergence rate as that of PPA.
%\end{Remark}

\subsection{Primal-dual algorithm via EG method}
We also develop a projection-based primal-dual algorithm via EG method to solve the problem \eqref{1}. Firstly, we compute the mid-point iteration $(x_{k+\frac{1}{2}},y_{k+\frac{1}{2}})$, i.e.,
\begin{subequations}\label{g1}
\begin{align}
&x_{k+\frac{1}{2}}=\mathcal{P}_{\mathcal{X}}\Big (x_{k}-\alpha \nabla_{x} f(x_{k},y_{k})\Big), \\
&y_{k+\frac{1}{2}}=\mathcal{P}_{\mathcal{Y}}\Big (y_{k}+\alpha \nabla_{y} f(x_{k},y_{k})\Big),
\end{align}
\end{subequations}
where $\alpha$ is the constant step-size that will be specified later. By using the mid-point $(x_{k+\frac{1}{2}},y_{k+\frac{1}{2}})$, we further compute the next iteration $(x_{k+1}, y_{k+1})$ as
\vspace{-6pt}
\begin{subequations}\label{g2}
\begin{align}
&x_{k+1}=\mathcal{P}_{\mathcal{X}}\Big(x_{k}-\alpha \nabla_{x} f(x_{k+\frac{1}{2}},y_{k+\frac{1}{2}})\Big), \\
&y_{k+1}=\mathcal{P}_{\mathcal{Y}}\Big(y_{k}+\alpha \nabla_{y} f(x_{k+\frac{1}{2}},y_{k+\frac{1}{2}})\Big).
\end{align}
\end{subequations}
According to the definitions of $z$ and $F(z)$ in \eqref{d2}, we can rewrite the algorithm \eqref{g1}-\eqref{g2} as
\vspace{-6pt}
\begin{subequations}\label{g3}
\begin{align}
& z_{k+\frac{1}{2}}=\mathcal{P}_{\Lambda}\Big(z_{k}-\alpha F(z_{k})\Big), \label{g3.1}\\
& z_{k+1}=\mathcal{P}_{\Lambda}\Big(z_{k}-\alpha F(z_{k+\frac{1}{2}})\Big). \label{g3.2}
\end{align}
\end{subequations}
It follows from \eqref{g3} that the crucial idea of the EG method is to find a mid-point $z_{k+\frac{1}{2}}$ by using the GDA method at the current point, and then obtain the next iteration by using the gradient $F(z_{k+\frac{1}{2}})$ at this mid-point. Compared with the GDA method in \cite{27}, the EG-based algorithm \eqref{g3} via adding the midpoint step can achieve exact convergence to a saddle point for general convex-concave functions.

\begin{Remark}
In contrast to the work of \cite{26}, the main differences of our proposed algorithms are two-fold. (i) We consider the constrained saddle-point problem while \cite{26} studied the unconstrained one. (ii) Our algorithms achieve exact convergence to a saddle point while the result of \cite{26} only converges into a bounded neighborhood of a saddle point.
\end{Remark}

%\begin{Remark}
%
%\end{Remark}
%This illustrates that the EG-based algorithm \eqref{g3} has double gradient computations and projection computations than those of ODGA-based algorithm \eqref{d2} and GDA method of \cite{27} at each iteration.
\subsection{Convergence analysis}
The convergence analyses of the proposed two primal-dual algorithms via OGDA and EG are provided. Firstly, we show the convergence result for the algorithm \eqref{d2} in the following theorem and its proof can be found in Appendix A.

\begin{Theorem}
Suppose that Assumptions 2-3 hold and the step-size $\alpha$ satisfies $0<\alpha< \frac{1}{2\kappa_{m}}$ with $\kappa_{m}=2\max(l_{xx},$ $l_{xy},l_{yx},l_{yy})$. Under the initial conditions $x_{0}=x_{-1}$ and $y_{0}=y_{-1}$, the developed OGDA-based algorithm \eqref{d2} guarantees that the iteration sequence $\{x_{k},y_{k}\}$ converges to a saddle point of problem \eqref{1}. Moreover, it holds that for any $T\ge 1$
\begin{align}\label{th1}
\vert f(\hat{x}_{T},\hat{y}_{T})-f(x^{*},y^{*}) \vert \le \frac{1}{2\alpha T}\Vert z_{0}-z^{*}\Vert^{2},
\end{align}
where $\hat{x}_{T}=\frac{1}{T}\sum^{T}_{k=1}x_{k}$ and $\hat{y}_{T}=\frac{1}{T}\sum^{T}_{k=1}y_{k}$.
\end{Theorem}

We next provide the convergence result of the algorithm \eqref{g3} with its proof given in Appendix B.

\begin{Theorem}
Suppose that Assumptions 2-3 hold and the step-size $\alpha$ satisfies $0<\alpha<\frac{1}{\kappa_{m}}$. Under the initial condition $z_{0}=z_{-1}$, the developed EG-based algorithm \eqref{g3} guarantees that the iteration sequence $\{x_{k},y_{k}\}$ converges to a saddle point of problem \eqref{1}. Furthermore, it holds that for any $T\ge 1$
\begin{align}\label{th2}
\vert f(\hat{x}_{T},\hat{y}_{T})-f(x^{*},y^{*}) \vert \le \frac{1}{2\alpha T}\Vert z_{0}-z^{*}\Vert^{2},
\end{align}
where $\hat{x}_{T}=\frac{1}{T}\sum^{T-1}_{k=0}x_{k+\frac{1}{2}}$ and $\hat{y}_{T}=\frac{1}{T}\sum^{T-1}_{k=0}y_{k+\frac{1}{2}}$.
\end{Theorem}

\begin{Remark}
It follows from Theorems 1-2 that the proposed OGDA-based algorithm \eqref{d2} and EG-based algorithm \eqref{g3} both achieve exact convergence to a saddle point rather than a bounded neighborhood of a saddle point shown in \cite{26}. Moreover, based on  \eqref{th1} and \eqref{th2} in Theorems 1-2, we have that the objective function $f(x,y)$ at the average iteration generated by these two algorithms converge to an optimal value with a sublinear rate $O(1/T)$.
\end{Remark}

\section{Unified Distributed algorithm via saddle-point dynamics}
Based on the primal-dual algorithms via OGDA and EG methods for constrained saddle-point problems, we develop unified distributed algorithm design and convergence analysis for solving the networked optimization problems \eqref{2} and \eqref{3}.

\subsection{Distributed constrained optimal consensus problem}
Note that the constrained optimal consensus problem \eqref{2} can be transformed into the constrained saddle-point problem \eqref{2a}. Based on the proposed OGDA-based algorithm \eqref{d1}, we develop a distributed primal-dual algorithm as
\begin{subequations}\label{d3}
\begin{align}
x^{k+1}_{i}&=\mathcal{P}_{\Omega_{i}}\Big (x^{k}_{i}-2\alpha \nabla f_{i}(x^{k}_{i})+\alpha \nabla f_{i}(x^{k-1}_{i})\\
&\quad-2\alpha \sum_{j\in \mathcal{N}_{i}} (x^{k}_{i}-x^{k}_{j}+v^{k}_{i}-v^{k}_{j}) \notag \\
&\quad+\alpha \sum_{j\in \mathcal{N}_{i}} (x^{k-1}_{i}-x^{k-1}_{j}+v^{k-1}_{i}-v^{k-1}_{j})\Big ), \notag \\
v^{k+1}_{i}&=v^{k}_{i}+2\alpha \sum_{j\in \mathcal{N}_{i}}(x^{k}_{i}-x^{k}_{j})-\alpha \sum_{j\in \mathcal{N}_{i}} (x^{k-1}_{i}-x^{k-1}_{j}).
\end{align}
\end{subequations}
Let $x_{k}=\text{col}(x^{k}_{i})^{N}_{i=1}, v_{k}=\text{col}(v^{k}_{i})^{N}_{i=1}, \nabla f(x_{k})=\text{col}(\nabla f_{i}(x^{k}_{i}))^{N}_{i=1}$, and $\mathcal{P}_{\Omega}(\cdot)=\text{col}(\mathcal{P}_{\Omega_{i}}(\cdot))^{N}_{i=1}$. From the definition of $L_{1}(x, v)$ in Section II, one has that $\nabla_{x} L_{1}(x,v)=\nabla f(x)+(L\otimes I_{m})(x+v)$ and $\nabla_{v} L_{1}(x,v)=(L\otimes I_{m})x$. Then, a compact form of \eqref{d3} can be obtained as
\begin{subequations}\label{d5}
\begin{align}
x_{k+1}&=\mathcal{P}_{\Omega}\Big (x_{k}-2\alpha \nabla_{x} L_{1}(\varpi_{k})+\alpha \nabla_{x} L_{1}(\varpi_{k-1}) \Big),\\
v_{k+1}&=v_{k}+2\alpha \nabla_{v} L_{1}(\varpi_{k})-\alpha \nabla_{v} L_{1}(\varpi_{k-1}),
\end{align}
\end{subequations}
where $\varpi_{k}=\text{col}(x_{k},v_{k})$. Define $\Phi(\varpi)=[ \nabla_{x} L_{1}(x,v);$ $-\nabla_{v} L_{1}(x,v)]=[\nabla f(x)+(L\otimes I_{m})(x+v);
-(L\otimes I_{m})x]$, and then algorithm \eqref{d5} can be arranged as $\varpi_{k+1}=\mathcal{P}_{\Theta_{1}}(\varpi_{k}-2\alpha \Phi(\varpi_{k})+\alpha \Phi(\varpi_{k-1}))$ with $\Theta_{1}=\Omega\times \mathbb{R}^{Nm}$. This illustrates that algorithm \eqref{d5} has the same structure as \eqref{d1}. Thus, the results of \eqref{d1} given in Theorem 1 can be easily extended to the case of \eqref{d5}. Under Assumption 3, one has that $\Phi(\varpi)$ is Lipschitz continuous, i.e., $\Vert \Phi(\varpi_{1})-\Phi(\varpi_{2})\Vert \le \kappa_{c} \Vert \varpi_{1}-\varpi_{2}\Vert$ for any $\varpi_{1}, \varpi_{2}$, where $\kappa_{c}$ is determined by Lipschitz constants of $\nabla f_{i}(x_{i}), i\in \mathcal{I}_{N}$ and the largest eigenvalue of $L$. Similar to the results of Theorem 1, we obtain the following corollary.

\begin{Corollary}
Suppose that Assumption 1 holds and the step-size $\alpha$ satisfies $0<\alpha<\frac{1}{2\kappa_{c}}$. The developed distributed algorithm \eqref{d3} guarantees that $x_{k}$ converges to an optimal solution of problem \eqref{2}. Moreover, for any $T\ge 1$, it holds that $\vert L_{1}(\hat{x}_{T},\hat{v}_{T})-L_{1}(x^{*},v^{*})\vert \le \frac{1}{2\alpha T}(\Vert x_{0}-x^{*}\Vert^{2}+\Vert v_{0}-v^{*}\Vert^{2})$, where $\hat{x}_{T}=\frac{1}{T}\sum^{T}_{k=1}x_{k}$ and $\hat{v}_{T}=\frac{1}{T}\sum^{T}_{k=1}v_{k}$.
\end{Corollary}
\vspace{-4pt}
\begin{Remark}
By applying the EG-based algorithm \eqref{g1}-\eqref{g2}, we develop another distributed primal-dual algorithm to solve the optimization problem \eqref{2}, which is composed of two steps.

\textbf{Step 1:} Calculate the mid-point iteration $(x^{k+\frac{1}{2}}_{i},v^{k+\frac{1}{2}}_{i})$.
\vspace{-2pt}
\begin{subequations}\label{hy1}
\begin{align}
x^{k+\frac{1}{2}}_{i}&=\mathcal{P}_{\Omega_{i}}\Big (x^{k}_{i}-\alpha \nabla f_{i}(x^{k}_{i})  \\
&\quad-\alpha \sum_{j\in \mathcal{N}_{i}} (x^{k}_{i}-x^{k}_{j}+v^{k}_{i}-v^{k}_{j})\Big ), \notag \\
v^{k+\frac{1}{2}}_{i}&=v^{k}_{i}+\alpha \sum_{j\in \mathcal{N}_{i}}(x^{k}_{i}-x^{k}_{j}).
\end{align}
\end{subequations}

\textbf{Step 2:} Calculate the next iteration $(x^{k+1}_{i},v^{k+1}_{i})$.
\vspace{-2pt}
\begin{subequations}\label{hy2}
\begin{align}
x^{k+1}_{i}&=\mathcal{P}_{\Omega_{i}}\Big (x^{k}_{i}-\alpha \nabla f_{i}(x^{k+\frac{1}{2}}_{i})  \\
&\quad-\alpha \sum_{j\in \mathcal{N}_{i}} (x^{k+\frac{1}{2}}_{i}-x^{k+\frac{1}{2}}_{j}+v^{k+\frac{1}{2}}_{i}-v^{k+\frac{1}{2}}_{j})\Big ), \notag \\
v^{k+1}_{i}&=v^{k}_{i}+\alpha \sum_{j\in \mathcal{N}_{i}}(x^{k+\frac{1}{2}}_{i}-x^{k+\frac{1}{2}}_{j}).
\end{align}
\end{subequations}
The proposed algorithm \eqref{hy1}-\eqref{hy2} obtains the same convergence results as Corollary 1, and its detailed proof can be derived from that of Theorem 2.
\end{Remark}

\begin{Remark}
%From the developed algorithm (), it seems that the neighbors¡¯ state $(x_{j}, v_{j})$ at the current step $k$ and last step $k-1$ both are required to be communicated, which leads to double communication than the results in []. Actually, at the current step $k$, the neighbors' variable $(x^{k-1}_{j}, v^{k-1}_{j})$
From the algorithm \eqref{d3}, it seems that the neighbors' states $(x_{j}, v_{j})$ at the current iteration and previous iteration are both transmitted, which leads to twice communication than those of \cite{7} and \cite{7a}. In fact, at the current iteration $k$, only $(x^{k}_{j}, v^{k}_{j})$ is required to be transmitted since $(x^{k-1}_{j}, v^{k-1}_{j})$ has been transmitted in the previous iteration. Thus, the communication requirement of the proposed algorithm \eqref{d3} is the same as those of \cite{7} and \cite{7a}.
\end{Remark}

\subsection{Distributed resource allocation problem}
Based on the OGDA-based algorithm \eqref{d1}, we propose a distributed algorithm to solve the optimization problem \eqref{3}
\begin{align*}
y^{k+1}_{i}&=\mathcal{P}_{\Omega_{i}}\Big ( y^{k}_{i}-2\alpha (\nabla h_{i}(y^{k}_{i})+W^{T}_{i}\lambda^{k}_{i})\\
&\quad+\alpha (\nabla h_{i}(y^{k-1}_{i})+W^{T}_{i}\lambda^{k-1}_{i})\Big), \tag{18a}\\
z^{k+1}_{i}&=z^{k}_{i}+2\alpha \sum_{j\in \mathcal{N}_{i}}(\lambda^{k}_{i}-\lambda^{k}_{j})\\
&\quad-\alpha  \sum_{j\in \mathcal{N}_{i}}(\lambda^{k-1}_{i}-\lambda^{k-1}_{j}), \tag{18b}
\end{align*}
\begin{subequations}
\begin{align*}
\lambda^{k+1}_{i}&=\lambda^{k}_{i}+2\alpha\Big (W_{i}y^{k}_{i}-d_{i}-\sum_{j\in \mathcal{N}_{i}}(z^{k}_{i}-z^{k}_{j}\\
&\quad+\lambda^{k}_{i}-\lambda^{k}_{j})\Big )-\alpha \Big (W_{i}y^{k-1}_{i}-d_{i} \\
&\quad-\sum_{j\in \mathcal{N}_{i}}(z^{k-1}_{i}-z^{k-1}_{j}+\lambda^{k-1}_{i}-\lambda^{k-1}_{j})\Big). \tag{18c}
\end{align*}
\end{subequations}
Let $y_{k}=\text{col}(y^{k}_{i})^{N}_{i=1}\in \mathbb{R}^{q}$, $z_{k}=\text{col}(z^{k}_{i})^{N}_{i=1}\in \mathbb{R}^{Nm}$, $\lambda_{k}=\text{col}(\lambda^{k}_{i})^{N}_{i=1}\in \mathbb{R}^{Nm}$, $\nabla h(y_{k})=\text{col}(\nabla h_{i}(y^{k}_{i}))^{N}_{i=1} \in \mathbb{R}^{q}$, $W=\text{diag}(W_{i})^{N}_{i=1}\in \mathbb{R}^{Nm\times q}$, and $d=\text{col}(d_{i})^{N}_{i=1}\in \mathbb{R}^{Nm}$. According to the definition of $L_{2}(y,z,\lambda)$ in Section II, we have that $\nabla_{y} L_{2}(y,z,\lambda)=\nabla h(y)+W^{T}\lambda$, $\nabla_{z} L_{2}(y,z,\lambda)=-(L\otimes I_{m})\lambda$, and $\nabla_{\lambda} L_{2}(y,z,\lambda)=Wy-d-(L\otimes I_{m})(z+\lambda)$. Then, a compact form of (18) is written as
\begin{subequations}\label{d8}
\begin{align}
y_{k+1}&=\mathcal{P}_{\Omega}\Big (y_{k}-2\alpha \nabla_{y} L_{2}(\xi_{k})+\alpha \nabla_{y} L_{2}(\xi_{k-1})\Big),\\
z_{k+1}&=z_{k}-2\alpha \nabla_{z} L_{2}(\xi_{k})+\alpha \nabla_{z} L_{2}(\xi_{k-1}),\\
\lambda_{k+1}&=\lambda_{k}+2\alpha \nabla_{\lambda} L_{2}(\xi_{k})-\alpha \nabla_{\lambda} L_{2}(\xi_{k-1}),
\end{align}
\end{subequations}
where $\xi_{k}=\text{col}(y_{k},z_{k},\lambda_{k})$. Define $\Psi(\xi)=[\nabla_{y}L_{2}(y,z,\lambda); $ $\nabla_{z} L_{2}(y,z,\lambda);-\nabla_{\lambda} L_{2}(y,z,\lambda)]=[\nabla h(y)+W^{T}\lambda;
-(L\otimes I_{m})\lambda; -(Wy-d-(L\otimes I_{m})(z+\lambda))]$, and \eqref{d8} is rewritten as $\xi_{k+1}=\mathcal{P}_{\Theta_{2}}(\xi_{k}-2\alpha \Psi(\xi_{k})+\alpha \Psi(\xi_{k-1}))$ with $\Theta_{2}=\Omega\times \mathbb{R}^{Nm}\times \mathbb{R}^{Nm}$, which has the same structure of \eqref{d1}. In addition, we obtain that $\Psi(\xi)$ is $\kappa_{s}$-Lipschitz continuous, where $\kappa_{s}$ is determined by Lipschitz constants of $\nabla h_{i}(y_{i}), i\in \mathcal{I}_{N}$ and largest eigenvalues of matrices $L$ and $W$.
\begin{Corollary}
Under Assumption 1 and the step-size satisfying $0<\alpha<\frac{1}{2\kappa_{s}}$, the distributed algorithm (18) guarantees that $y_{k}$ converges to an optimal solution of the problem \eqref{3}. Moreover, $\vert L_{2}(\hat{y}_{T},\hat{z}_{T}, \hat{\lambda}_{T})-L_{2}(x^{*},z^{*},\lambda^{*})\vert \le \frac{1}{2\alpha T}(\Vert y_{0}-y^{*}\Vert^{2}+\Vert z_{0}-z^{*}\Vert^{2}+\Vert \lambda_{0}-\lambda^{*}\Vert^{2})$ holds for any $T\ge 1$,, where $\hat{y}_{T}=\frac{1}{T}\sum^{T}_{k=1}y_{k}$, $\hat{z}_{T}=\frac{1}{T}\sum^{T}_{k=1}z_{k}$ and $\hat{\lambda}_{T}=\frac{1}{T}\sum^{T}_{k=1}\lambda_{k}$.
\end{Corollary}
\begin{Remark}
Based on the EG-based algorithm \eqref{g1}-\eqref{g2}, another distributed primal-dual algorithm is developed to solve the optimization problem \eqref{3}, which is formulated as

\textbf{Step 1:} Calculate the mid-point $(y^{k+\frac{1}{2}}_{i},z^{k+\frac{1}{2}}_{i}, \lambda^{k+\frac{1}{2}}_{i})$.
\begin{subequations}\label{hy3}
\begin{align}
y^{k+\frac{1}{2}}_{i}&=\mathcal{P}_{\Omega_{i}}\Big ( y^{k}_{i}-\alpha (\nabla h_{i}(y^{k}_{i})+W^{T}_{i}\lambda^{k}_{i})\Big), \\
z^{k+\frac{1}{2}}_{i}&=z^{k}_{i}+\alpha \sum_{j\in \mathcal{N}_{i}}(\lambda^{k}_{i}-\lambda^{k}_{j}),\\
\lambda^{k+\frac{1}{2}}_{i}&=\lambda^{k}_{i}+\alpha \Big (W_{i}y^{k}_{i}-d_{i} \notag \\
&\quad-\sum_{j\in \mathcal{N}_{i}}(z^{k}_{i}-z^{k}_{j}+\lambda^{k}_{i}-\lambda^{k}_{j}) \Big).
\end{align}
\end{subequations}

\textbf{Step 2:} Calculate the next iteration $(y^{k+1}_{i},z^{k+1}_{i}, \lambda^{k+1}_{i})$.
\begin{subequations}\label{hy4}
\begin{align}
y^{k+1}_{i}&=\mathcal{P}_{\Omega_{i}}\Big ( y^{k+\frac{1}{2}}_{i}-\alpha (\nabla h_{i}(y^{k+\frac{1}{2}}_{i})+W^{T}_{i}\lambda^{k+\frac{1}{2}}_{i})\Big),  \\
z^{k+1}_{i}&=z^{k}_{i}+\alpha \sum_{j\in \mathcal{N}_{i}}(\lambda^{k+\frac{1}{2}}_{i}-\lambda^{k+\frac{1}{2}}_{j}), \\
\lambda^{k+1}_{i}&=\lambda^{k}_{i}+\alpha \Big (W_{i}y^{k+\frac{1}{2}}_{i}-d_{i} \notag\\
&\quad-\sum_{j\in \mathcal{N}_{i}}(z^{k+\frac{1}{2}}_{i}-z^{k+\frac{1}{2}}_{j}+\lambda^{k+\frac{1}{2}}_{i}-\lambda^{k+\frac{1}{2}}_{j}) \Big ).
\end{align}
\end{subequations}
Actually, the algorithm \eqref{hy3}-\eqref{hy4} has the same formulation as that in \cite{17}. However, only asymptotic convergence was proven in \cite{17} and its convergence rate analysis was not given. Based on the result of Theorems 2, we easily prove that the algorithm \eqref{hy3}-\eqref{hy4} achieves exact convergence to an optimal solution with $O(1/k)$ convergence rate.
\end{Remark}
%Compared with the algorithms in \cite{14,15,16} with the diminishing step-sizes, the proposed distributed algorithms (18) and \eqref{hy3}-\eqref{hy4} both use the constant step-size and therefore the convergence result outperforms those of \cite{14,15,16}.

\begin{Remark}
Although the traditional
centralized optimization method (e.g., ADMM-based algorithm in \cite{31}) also can solve these two networked optimization problems, it requires massive communication and large bandwidth for the central node. In contrast, the developed distributed algorithm via local information interaction can overcome the issues of the centralized method and therefore can be applied to solve a large-scale networked optimization problem. In addition, unlike the distributed algorithms in \cite{4a} and \cite{14a} that require solving a sub-optimization problem at each iteration, the developed algorithms are easier to be implemented without solving the sub-optimization problem.
%In particular, the algorithm of \cite{13} only converges to a suboptimal solution, and the convergence rates of \cite{14,15,16} are slow since the diminishing step-size is used. The proposed distributed algorithms \eqref{d6} and \eqref{hy3}-\eqref{hy4} both use the constant step-size and therefore the convergence result outperforms those of \cite{14,15,16}. In
\end{Remark}
%
%\begin{Remark}
%
%\end{Remark}
\section{Numerical simulation}
In this section, we provide some numerical simulation examples for solving networked optimization problems \eqref{3} and general constrained saddle-point problem \eqref{1} to demonstrate the effectiveness of the proposed algorithms.

\emph{Example 1:} We first verify the proposed OGDA-based algorithms \eqref{d2} and EG-based algorithm \eqref{g3} by solving the following constrained saddle-point problem
\begin{align}\label{hy5}
\min_{x\in \mathcal{X}} \max_{y\in \mathcal{Y}} f(x,y)=x^{T}By,
\end{align}
where $B \in \mathbb{R}^{10\times 10}$ is a random matrix and its element is generated from a uniform distribution on $[0,5]$, the constrained sets $\mathcal{X}$ and $\mathcal{Y}$ are set to be $\mathcal{X}=[-5,5]^{10}$ and $\mathcal{Y}=[-2,2]^{10}$. Then, we obtain that $(x^{*}, y^{*})=(0_{10}, 0_{10})$ is a saddle point of problem \eqref{hy5} and the optimal value is $f(x^{*},y^{*})=0$.

We carry out the OGDA-based algorithm \eqref{d2} and EG-based algorithm \eqref{g3} under the same initial values $x_{0}=10\times\bm 1_{10}$, $y_{0}=10\times \bm 1_{10}$, and the chosen step-size $\alpha=0.01$. In addition, the GDA algorithm of \cite{27} is also implemented as a comparison. Fig. 1 shows that the convergence results of the objective error $\vert f(x_{k},y_{k})-f(x^{*},y^{*})\vert$ under the OGDA-based algorithm \eqref{d2}, EG-based algorithm \eqref{g3} and GDA algorithm of \cite{27}. It is shown that the developed OGDA algorithm \eqref{d2} and EG algorithm \eqref{g3} both guarantee that the iteration $(x_{k},y_{k})$ converges to the saddle point $(0_{10},0_{10})$ while the GDA algorithm of \cite{27} does not converge.

%\begin{figure}[!ht]
%\centering
%\includegraphics[width=0.45\textwidth, clip=true]{saddle_single1}
%\caption{$\vert f(x_{k},y_{k})-f(x^{*},y^{*})\vert$ under the OGDA-based algorithm \eqref{d2}, EG-based algorithm \eqref{g3}, and GDA algorithm of \cite{27}.}
%\end{figure}

\begin{figure}[!ht]
\centering
\includegraphics[width=0.5\textwidth, clip=true]{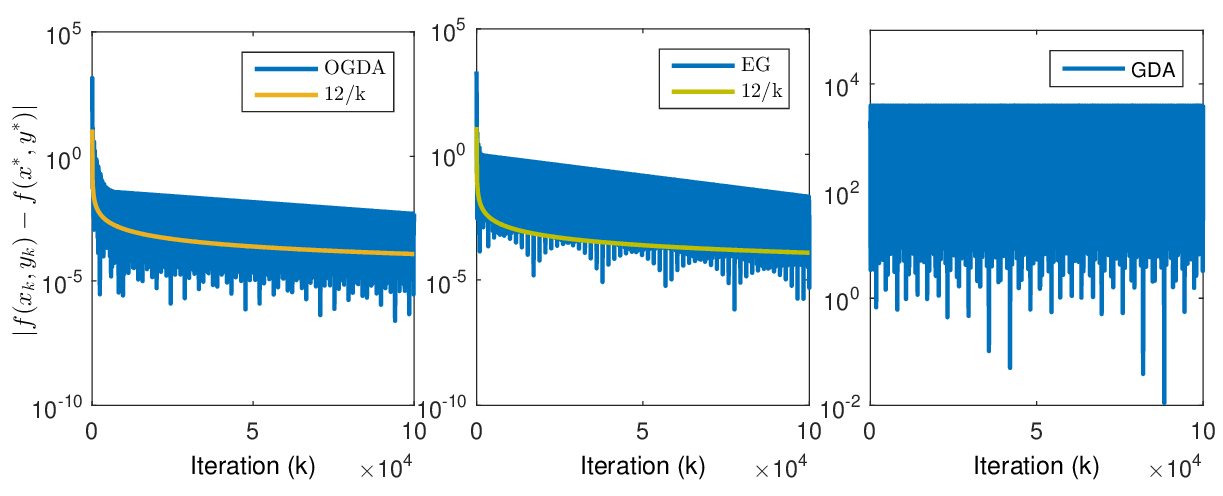}
\caption{$\vert f(x_{k},y_{k})-f(x^{*},y^{*})\vert$ under the OGDA-based algorithm \eqref{d2}, EG-based algorithm \eqref{g3} and GDA algorithm of \cite{27}.}
\end{figure}

\emph{Example 2:} We next demonstrate the distributed ODGA-based algorithm (18) and EG-based algorithm \eqref{hy3}-\eqref{hy4} to solve the resource allocation problem \eqref{3}. Consider a network of $N=20$ agents and its topology is described by a ring graph. Each local objective function is $h_{i}(y_{i})=a_{i}y_{i}+b_{i}\text{log}(1+e^{c_{i}y_{i}})$, and local set constraint is $\Omega_{i}=[-1,1], i\in \mathcal{I}_{N}$. The datums in the function $h_{i}(y_{i})$ and coupled equation constraint $\sum^{N}_{i=1}(W_{i}y_{i}-d_{i})=0$  are randomly generated from $a_{i}\in [-5,5], b_{i}\in [0,2], c_{i}\in [0,1], W_{i}\in [-1,1]$ and $d_{i} \in [-2,2]$.

\begin{figure}[!ht]
\centering
\includegraphics[width=0.5\textwidth, clip=true]{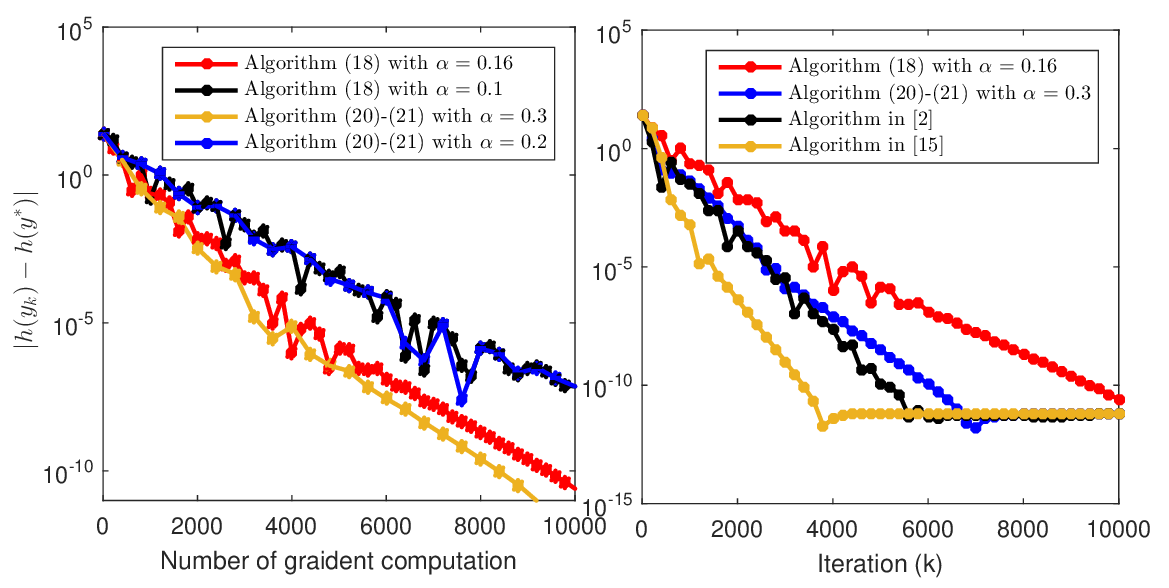}
\caption{The objective error $\vert h(y_{k})-h(y^{*})\vert$ under the distributed ODGA-based algorithm (18) and EG-based algorithm \eqref{hy3}-\eqref{hy4}}
\end{figure}

The distributed ODGA-based algorithm (18) and EG-based algorithm \eqref{hy3}-\eqref{hy4} are implemented by choosing different step-sizes $\alpha$. The left subfigure of Fig. 2 describes the objective error $\vert h(y_{k})-h(y^{*})\vert$ under these two algorithms with respect to the number of gradient computation, which shows that the developed two algorithms both guarantee $\vert h(y_{k})-h(y^{*})\vert$ converging to zero. The right subfigure of Fig. 2 provides the performance comparison between the developed two algorithms and the algorithms in \cite{4a} and \cite{14a}. It is shown that the algorithm in \cite{4a} enjoys better convergence performance than the other algorithms, and the developed EG-based algorithm has similar convergence performance as that in \cite{14a}.

% From It is illustrated that the convergence result of the ODGA-based algorithm \eqref{d6} outperforms that of EG-based algorithm \eqref{hy3}-\eqref{hy4} in terms of the number of gradient computations.

\section{Conclusion}
This paper develops a unified distributed method for solving two classes of networked optimization problems with non-identical set constraints. We first establish the relationship between two networked optimization problems and constrained saddle-point problems, and then propose two projection-based primal-dual algorithms via OGDA and EG methods. Subsequently, we develop unified distributed algorithms via saddle-point dynamics to solve these two networked optimization problems. The final examples demonstrates the effectiveness of the developed algorithms.

%Further works include considering a class of constrained non-convex and non-concave saddle-point problem and its application to network optimization.

\appendix
\renewcommand{\thesection}{Appendix A}
Before presenting the proofs of Theorems 1-2, some preliminary results are provided.

\begin{Lemma}[Lemma 4, \cite{26}]
Let $F(\cdot)$ be defined in \eqref{d2}. Under Assumption 1, the following results hold
\begin{enumerate}
  \item [\text{(i)}] $F(\cdot)$ is a monotone operator, i.e., $(F(z_{1})-F(z_{2}))^\mathrm{T}(z_{1}-z_{2})\ge 0$ for any $z_{1}, z_{2}\in \mathcal{X}\times \mathcal{Y}$.
  \item [\text{(ii)}] $F(\cdot)$ is Lipschitz continuous, i.e., $\Vert F(z_{1})-F(z_{2})\Vert \le \kappa_{m} \Vert z_{1}-z_{2}\Vert$ holds for any $z_{1}, z_{2}\in \mathcal{X}\times \mathcal{Y}$, where $\kappa_{m}=2\max(l_{xx},l_{xy},l_{yx},l_{yy})$.
\end{enumerate}
\end{Lemma}

\begin{Lemma}[Proposition 7, \cite{26}]
Let $\{z_{k}\}$ be the iteration sequence generated by the following update
\begin{align*}
z_{k+1}=z_{k}-\alpha F(z_{k+1})+\varepsilon_{k},
\end{align*}
where $F:\mathbb{R}^{m+n}\to \mathbb{R}^{m+n}$ is a continuous function, $\alpha$ is a positive constant, and $\varepsilon_{k}\in \mathbb{R}^{m+n}$ is an arbitrary vector. For any $z\in \mathbb{R}^{m+n}$ and $k\ge 1$, it holds that
\begin{align}\label{p0}
&\Vert z_{k+1}-z\Vert^{2}=\Vert z_{k}-z\Vert^{2}-2\alpha (z_{k+1}-z)^\mathrm{T}F(z_{k+1})  \\
&~~~~~~~~~~~-\Vert z_{k+1}-z_{k}\Vert^{2}+2\varepsilon^\mathrm{T}_{k}(z_{k+1}-z). \notag
\end{align}
\end{Lemma}

\begin{Lemma}[Proposition 5, \cite{26}]
Define $\hat{x}_{T}=\frac{1}{T}\sum^{T}_{k=1}x_{k}$ and $\hat{y}_{T}=\frac{1}{T}\sum^{T}_{k=1}y_{k}$. Under Assumption 1, it follows that
\begin{align}\label{lm1}
&f(\hat{x}_{T},y^{*})-f(x^{*},\hat{y}_{T}) \notag \\
&~~~~~~\le \frac{1}{T} \sum^{T-1}_{k=0} (z_{k+1}-z^{*})^{\mathrm{T}}F(z_{k+1}).
\end{align}
\end{Lemma}

\subsection{Proof of Theorem 1}
Define $Y_{k+1}=\mathcal{P}_{\Lambda}(\Upsilon_{k})-\Upsilon_{k}$, where $\Upsilon_{k}=z_{k}-\alpha F(z_{k})-\alpha(F(z_{k})-F(z_{k-1}))$. It then follows from \eqref{d2} that
\begin{align*}
z_{k+1}=z_{k}-\alpha F(z_{k})-\alpha(F(z_{k})-F(z_{k-1}))+Y_{k+1},
\end{align*}
which can be rewritten as $z_{k+1}=z_{k}-\alpha F(z_{k+1})+\chi_{k}$ with $\chi_{k}=\alpha \{(F(z_{k+1})-F(z_{k})-(F(z_{k})-F(z_{k-1}))\}+Y_{k+1}$. By using \eqref{p0} of Lemma 2, we obtain that
\begin{align}\label{p3}
&(z_{k+1}-z)^\mathrm{T}F(z_{k+1}) \notag \\
&=\frac{1}{2\alpha}\Vert z_{k}-z\Vert^{2}-\frac{1}{2\alpha}\Vert z_{k+1}-z\Vert^{2} \notag\\
&\quad-\frac{1}{2\alpha}\Vert z_{k+1}-z_{k}\Vert^{2}+\frac{1}{\alpha}\chi^\mathrm{T}_{k}(z_{k+1}-z) \notag\\
&=\frac{1}{2\alpha}\Vert z_{k}-z\Vert^{2}-\frac{1}{2\alpha}\Vert z_{k+1}-z\Vert^{2}-\frac{1}{2\alpha}\Vert z_{k+1}-z_{k}\Vert^{2} \notag\\
&\quad+\frac{1}{\alpha}Y^\mathrm{T}_{k+1}(z_{k+1}-z)+(z_{k+1}-z)^\mathrm{T}(F(z_{k+1})-F(z_{k})) \notag\\
&\quad-(z_{k}-z)^{\mathrm{T}}(F(z_{k})-F(z_{k-1})) \notag \\
&\quad-(z_{k+1}-z_{k})^{\mathrm{T}}(F(z_{k})-F(z_{k-1})).
\end{align}
According to the Lipschitz continuity of $F(z)$ in Lemma 1 and Young's inequality, we have that $-(z_{k+1}-z_{k})^{\mathrm{T}}(F(z_{k})-F(z_{k-1}))\le \frac{\kappa_{m}}{2} \Vert z_{k}-z_{k+1}\Vert^{2}+\frac{\kappa_{m}}{2} \Vert z_{k}-z_{k-1}\Vert^{2}$. Then, \eqref{p3} can be simplified as
\begin{align}\label{p6}
&(z_{k+1}-z)^\mathrm{T}F(z_{k+1}) \notag\\
&\le \frac{1}{2\alpha}\Vert z_{k}-z\Vert^{2}-\frac{1}{2\alpha}\Vert z_{k+1}-z\Vert^{2}-\eta \Vert z_{k+1}-z_{k}\Vert^{2}  \\
&\quad -\frac{\kappa_{m}}{2} \Vert z_{k+1}-z_{k}\Vert^{2}+\frac{\kappa_{m}}{2} \Vert z_{k}-z_{k-1}\Vert^{2} \notag\\
&\quad+(z_{k+1}-z)^\mathrm{T}(F(z_{k+1})-F(z_{k}))\notag\\
&\quad-(z_{k}-z)^\mathrm{T}(F(z_{k})-F(z_{k-1}))+\frac{1}{\alpha}Y^\mathrm{T}_{k+1}(z_{k+1}-z), \notag
\end{align}
where $\eta=\frac{1}{2\alpha}-\kappa_{m}>0$ if $\alpha<\frac{1}{2\kappa_{m}}$ is chosen.

Let $z^{*}=\text{col}(x^{*},y^{*})$ be a saddle point of problem \eqref{1}. Since $z_{k}\in \mathcal{X}\times \mathcal{Y}$ for all $k\ge 1$, it then follows from \eqref{1a} that $(z_{k}-z^{*})^\mathrm{T}F(z^{*})\ge 0, \forall k\ge 1$. According to the monotone property of $F(\cdot)$ given in Lemma 1, one has that $(z_{k}-z^{*})^\mathrm{T}(F(z_{k})-F(z^{*}))\ge 0$, which further implies that
\begin{align}\label{p7a}
(z_{k}-z^{*})^\mathrm{T}F(z_{k})\ge 0, \forall k\ge 1.
\end{align}
In addition, according to the definition of $Y_{k+1}$, one has that $\mathcal{P}_{\Lambda}(z_{k+1}-Y_{k+1})=z_{k+1}$. This implies that $-Y_{k+1}\in \mathcal{N}_{\Lambda}(z_{k+1})$, where $\mathcal{N}_{\Lambda}(z_{k+1})$ is the normal cone of the set $\Lambda$ at $z_{k+1}$. Since $z^{*}\in \Lambda$, one can further derive that
\begin{align}\label{p7}
(z_{k+1}-z^{*})^\mathrm{T}Y_{k+1}\le 0.
\end{align}
Setting $z=z^{*}$ of \eqref{p6}, and combining \eqref{p7a}-\eqref{p7}, one has that
\begin{align}\label{p8}
0& \le (z_{k+1}-z^{*})^\mathrm{T}F(z_{k+1})\le \frac{1}{2\alpha}\Vert z_{k}-z^{*}\Vert^{2}\\
&\quad-\frac{1}{2\alpha}\Vert z_{k+1}-z^{*}\Vert^{2}-\eta \Vert z_{k+1}-z_{k}\Vert^{2}-\frac{\kappa_{m}}{2} \Vert z_{k+1}-z_{k}\Vert^{2}\notag\\
&\quad+\frac{\kappa_{m}}{2} \Vert z_{k}-z_{k-1}\Vert^{2}+(z_{k+1}-z^{*})^\mathrm{T}(F(z_{k+1})-F(z_{k}))\notag\\
&\quad-(z_{k}-z^{*})^\mathrm{T}(F(z_{k})-F(z_{k-1})). \notag
\end{align}
Summing \eqref{p8} over $k$ from $0$ to $t$, we obtain that
\begin{align}\label{p9}
&\eta \sum^{t}_{k=0} \Vert z_{k+1}-z_{k}\Vert^{2} \le \frac{1}{2\alpha}\Vert z_{0}-z^{*}\Vert^{2}-\frac{1}{2\alpha}\Vert z_{t+1}-z^{*}\Vert^{2} \notag \\
& \quad -\frac{\kappa_{m}}{2} \Vert z_{t+1}-z_{t}\Vert^{2}+\frac{\kappa_{m}}{2} \Vert z_{0}-z_{-1}\Vert^{2}+(z_{t+1}-z^{*})^\mathrm{T} \notag \\
& \quad (F(z_{t+1})-F(z_{t}))-(z_{0}-z^{*})^\mathrm{T}(F(z_{0})-F(z_{-1})) \notag \\
&\le \frac{1}{2\alpha}\Vert z_{0}-z^{*}\Vert^{2}-(\frac{1}{2\alpha}-\frac{\kappa_{m}}{2})\Vert z_{t+1}-z^{*}\Vert^{2},
\end{align}
where the last second inequality is obtained by using the initial condition $z_{0}=z_{-1}$ and $(z_{t+1}-z^{*})^{T}(F(z_{t+1})-F(z_{t}))\le \frac{\kappa_{m}}{2}\Vert z_{t+1}-z^{*} \Vert^{2}+\frac{\kappa_{m}}{2}\Vert z_{t+1}-z_{t} \Vert^{2}$. Letting $t\to \infty$ and under $\alpha<\frac{1}{2\kappa_{m}}$, it follows from \eqref{p9} that
\begin{align*}
\sum^{\infty}_{k=0} \Vert z_{k+1}-z_{k}\Vert^{2}\le \frac{1}{2\alpha\eta}\Vert z_{0}-z^{*}\Vert^{2}<\infty.
\end{align*}
Consequently, we obtain that $\lim_{k\to \infty}(z_{k+1}-z_{k})=0$. In addition, it follows from \eqref{p9} that $(\frac{1}{2\alpha}-\frac{\kappa_{m}}{2})\Vert z_{t+1}-z^{*}\Vert^{2}\le \frac{1}{2\alpha}\Vert z_{0}-z^{*}\Vert^{2}$ holds for any $t\ge 0$. This implies that $z_{k}$ is bounded for $\forall k\in \mathbb{N}$. Then, we obtain that $z_{k}$ has the subsequence $\{z_{n_{k}}\}$ that converges to some limit point $z^{\infty}$, i.e., $\lim_{k\to \infty} z_{n_{k}}=z^{\infty}=\text{col}(x^{\infty},\lambda^{\infty})$. Moreover, from \eqref{d2}, we derive that
\begin{align*}
&z^{\infty}=\mathcal{P}_{\Lambda}(z^{\infty}-\alpha F(z^{\infty}))=0.
\end{align*}
This implies that $-\alpha F(z^{\infty}) \in \mathcal{N}_{\Lambda}(z^{\infty})$ and then we obtain that $(z-z^{\infty})^{T}F(z^{\infty})\ge 0$ holds for $z\in \mathcal{X}\times \mathcal{Y}$. It follows from \eqref{1a} that $z^{\infty}$ is a saddle point of problem \eqref{1}.

To this end, we have shown that $\{z_{k}\}$ has a convergence subsequence $\{z_{n_{k}}\}$. We next prove the convergence of the original sequence $\{z_{k}\}$. From \eqref{p8}, one has that $\frac{1}{2\alpha}\Vert z_{k+1}-z^{*}\Vert^{2}+\frac{\kappa_{m}}{2} \Vert z_{k+1}-z_{k}\Vert^{2}-(z_{k+1}-z^{*})^\mathrm{T}(F(z_{k+1})-F(z_{k}))\le \frac{1}{2\alpha}\Vert z_{k}-z^{*}\Vert^{2}+\frac{\kappa_{m}}{2} \Vert z_{k}-z_{k-1}\Vert^{2}-(z_{k}-z^{*})^\mathrm{T}(F(z_{k})-F(z_{k-1}))-\eta \Vert z_{k+1}-z_{k}\Vert^{2}$. Define $\Delta_{k}=\frac{1}{2\alpha}\Vert z_{k}-z^{*}\Vert^{2}-\frac{\kappa_{m}}{2} \Vert z_{k}-z_{k-1}\Vert^{2}+(z_{k}-z^{*})^\mathrm{T}(F(z_{k})-F(z_{k-1}))$ and one has that $\Delta_{k}\ge (\frac{1}{2\alpha}-\frac{\kappa_{m}}{2})\Vert z_{k}-z^{*}\Vert^{2}\ge 0$. It then follows that
\begin{align*}
\Delta_{k+1}\le \Delta_{k}-\eta \Vert z_{k+1}-z_{k}\Vert^{2}.
\end{align*}
According to the monotonicity and boundedness of $\Delta_{k}$, we have that $\Delta_{k}$ is convergent. Based on the fact that $\lim_{k \to \infty} (z_{k+1}-z_{k})=0$, one has that $\Vert z_{k}-z^{*}\Vert$ is convergent. By setting $z^{*}=z^{\infty}$, we have that $\lim_{k\to \infty} z_{n_{k}}=z^{\infty}=z^{*}$. Based on $\lim_{k\to \infty} z_{n_{k}}=z^{*}$ and $\lim_{k\to \infty}(z_{k+1}-z_{k})=0$, we obtain that $\lim_{k\to \infty} \Delta_{k}=0$. Under the fact that $\Delta_{k}\ge (\frac{1}{2\alpha}-\frac{\kappa_{m}}{2})\Vert z_{k}-z^{*}\Vert^{2}\ge 0$, we obtain that $\lim_{k\to \infty} \Vert z_{k}-z^{*}\Vert^{2}=0$. Thus, we have shown that the sequence $\{z_{k}\}$ converges to a saddle point of problem \eqref{1}.

We next analyze the convergence rate of algorithm \eqref{d2}. From \eqref{lm1} in Lemma 3, we obtain that
\begin{align*}
&  f(\hat{x}_{T},y^{*})-f(x^{*},\hat{y}_{T})\le \frac{1}{T} \sum^{T-1}_{k=0} (z_{k+1}-z^{*})^\mathrm{T}F(z_{k+1}) \notag \\
& \le \frac{1}{T}\Big (\frac{1}{2\alpha}\Vert z_{0}-z^{*}\Vert^{2}-\frac{1}{2\alpha}\Vert z_{T}-z^{*}\Vert^{2}-\frac{\kappa_{m}}{2} \Vert z_{T} \notag \\
&\quad-z_{T-1}\Vert^{2}+(z_{T}-z^{*})^\mathrm{T}(F(z_{T})-F(z_{T-1}))\Big ) \notag\\
& \le \frac{1}{T}\Big (\frac{1}{2\alpha}\Vert z_{0}-z^{*}\Vert^{2}-(\frac{1}{2\alpha}-\frac{\kappa_{m}}{2})\Vert z_{t+1}-z^{*}\Vert^{2} \Big ) \notag \\
& \le \frac{1}{2\alpha T}\Vert z_{0}-z^{*}\Vert^{2}.
\end{align*}
where the second inequality is derived from \eqref{p8} and the third inequality is obtained by using \eqref{p9}. Note that $f(\hat{x}_{T},y^{*})-f(x^{*},\hat{y}_{T})=f(\hat{x}_{T},y^{*})-f(x^{*},y^{*})+f(x^{*}, y^{*})-f(x^{*},\hat{y}_{T})\le \frac{1}{2\alpha T}\Vert z_{0}-z^{*}\Vert^{2}$. Since $(\hat{x}_{T},\hat{y}_{T}) \in \mathcal{X}\times \mathcal{Y}$ and $f(x,y)$ is a convex and concave function on $\mathcal{X}\times \mathcal{Y}$, we obtain that $0 \le f(\hat{x}_{T},y^{*})-f(x^{*},y^{*})\le  \frac{1}{2\alpha T}\Vert z_{0}-z^{*}\Vert^{2}$ and $0 \le f(x^{*}, y^{*})-f(x^{*},\hat{y}_{T})\le \frac{1}{2\alpha T}\Vert z_{0}-z^{*}\Vert^{2}$.
%\begin{align*}
%&0 \le f(\hat{x}_{N},\lambda^{*})-f(x^{*},\lambda^{*})\le  \frac{1}{2\alpha N}\Vert z_{0}-z^{*}\Vert^{2}, \\
%&0 \le f(x^{*}, \lambda^{*})-f(x^{*},\hat{\lambda}_{N})\le \frac{1}{2\alpha N}\Vert z_{0}-z^{*}\Vert^{2}.
%\end{align*}
In addition, since $f(\hat{x}_{T},\hat{y}_{T})\le f(\hat{x}_{T},y^{*})$ and $f(\hat{x}_{T},\hat{y}_{T})\ge f(x^{*},\hat{y}_{T})$, we further derive that $f(\hat{x}_{T},\hat{y}_{T})-f(x^{*},y^{*}) \le f(\hat{x}_{T},y^{*})-f(x^{*},y^{*})\le \frac{1}{2\alpha T}\Vert z_{0}-z^{*}\Vert^{2}$ and $f(x^{*}, y^{*})-f(\hat{x}_{T},\hat{y}_{T})\le f(x^{*}, y^{*})-f(x^{*},\hat{y}_{T})\le  \frac{1}{2\alpha T}\Vert z_{0}-z^{*}\Vert^{2}$. Thus, we obtain that $\vert f(\hat{x}_{T},\hat{y}_{T})-f(x^{*},y^{*}) \vert \le \frac{1}{2\alpha T}\Vert z_{0}-z^{*}\Vert^{2}$.
\vspace{-0.6pt}
\subsection{Proof of Theorem 2}
Let $m_{k+1}=z_{k+\frac{1}{2}}-(z_{k}-\alpha F(z_{k}))$ and we obtain that $z_{k+\frac{1}{2}}=z_{k}-\alpha F(z_{k})+m_{k+1}$. It then follows from \eqref{g3.1} that $\mathcal{P}_{\Lambda}(z_{k+\frac{1}{2}}-m_{k+1})=z_{k+\frac{1}{2}}$, which implies $-m_{k+1}\in \mathcal{N}_{\Lambda}(z_{k+\frac{1}{2}})$. Since $z_{k+1}\in \Lambda$, one has that $(z_{k+\frac{1}{2}}-z_{k+1})^{T}m_{k+1}\le 0$. In addition, define $n_{k+1}=z_{k+1}-(z_{k}-\alpha F(z_{k+\frac{1}{2}}))$ and one can derive that
\vspace{-0.6pt}
\begin{align}\label{p10}
z_{k+1}=z_{k}-\alpha F(z_{k+\frac{1}{2}})+n_{k+1}.
\end{align}
Note from \eqref{g3.2} that $\mathcal{P}_{\Lambda}(z_{k+1}-n_{k+1})=z_{k+1}$, which infers $-n_{k+1}\in \mathcal{N}_{\Lambda}(z_{k+1})$. It then follows that $(z_{k+1}-z^{*})^{T}n_{k+1}\le 0$. Also, the above equation \eqref{p10} can be rewritten as
\vspace{-0.6pt}
\begin{align}\label{p11}
z_{k+1}=z_{k}-\alpha F(z_{k+1})+\psi_{k+1},
\end{align}
where $\psi_{k+1}=\alpha F(z_{k+1})-\alpha F(z_{k+\frac{1}{2}})+n_{k+1}$.

By setting $z=z^{*}$ of \eqref{p0} in Lemma 2, it follows from \eqref{p11} that
\begin{align*}
&\Vert z_{k+1}-z^{*}\Vert^{2}=\Vert z_{k}-z^{*}\Vert^{2}-2\alpha (z_{k+1}-z^{*})^{T}F(z_{k+1}) \notag\\
&~~~~-\Vert z_{k+1}-z_{k}\Vert^{2}+2\psi^{T}_{k+1}(z_{k+1}-z^{*})
\end{align*}
\begin{align}\label{p12}
&\le \Vert z_{k}-z^{*}\Vert^{2}-\Vert z_{k+1}-z_{k}\Vert^{2}-2\alpha (z_{k+1}-z^{*})^{T}F(z_{k+\frac{1}{2}}) \notag \\
&\le \Vert z_{k}-z^{*}\Vert^{2}-\Vert z_{k+1}-z_{k+\frac{1}{2}}\Vert^{2}-\Vert z_{k+\frac{1}{2}}-z_{k}\Vert^{2} \notag\\
&\quad-2(z_{k+1}-z_{k+\frac{1}{2}})^{T}(-\alpha F(z_{k})+m_{k+1}) \notag\\
&\quad -2\alpha (z_{k+1}-z_{k+\frac{1}{2}})^{T}F(z_{k+\frac{1}{2}})-2\alpha (z_{k+\frac{1}{2}}-z^{*})^{T}F(z_{k+\frac{1}{2}}) \notag\\
& \le \Vert z_{k}-z^{*}\Vert^{2}-\Vert z_{k+1}-z_{k+\frac{1}{2}}\Vert^{2}-\Vert z_{k+\frac{1}{2}}-z_{k}\Vert^{2} \notag\\
&\quad-2\alpha(z_{k+1}-z_{k+\frac{1}{2}})^{T}(F(z_{k+\frac{1}{2}})-F(z_{k})) \notag\\
&\quad-2\alpha (z_{k+\frac{1}{2}}-z^{*})^{T}F(z_{k+\frac{1}{2}}),
\end{align}
where the first inequality is obtained by using $(z_{k+1}-z^{*})^{T}n_{k+1}\le 0$, and the second inequality is derived with $\Vert a-b \Vert^{2}=\Vert a-c \Vert^{2}+\Vert b-c \Vert^{2}+2(a-c)^{T}(c-b)$ for any vector $a, b, c$ and $z_{k+\frac{1}{2}}=z_{k}-\alpha F(z_{k})+m_{k+1}$, and the last inequality is obtained by using $(z_{k+\frac{1}{2}}-z_{k+1})^{T}m_{k+1}\le 0$. Note that $-2\alpha(z_{k+1}-z_{k+\frac{1}{2}})^{T}(F(z_{k+\frac{1}{2}})-F(z_{k}))\le \alpha^{2}\kappa^{2}_{m} \Vert z_{k+1}-z_{k+\frac{1}{2}}\Vert^{2}+\Vert z_{k+\frac{1}{2}}-z_{k}\Vert^{2}$, and it then follows from \eqref{p12} that
\begin{align}\label{p13}
&(z_{k+\frac{1}{2}}-z^{*})^{T} F(z_{k+\frac{1}{2}})\le \frac{1}{2\alpha}\Vert z_{k}-z^{*}\Vert^{2}\\
&\quad~~~~~~~~-\frac{1}{2\alpha}\Vert  z_{k+1}-z^{*}\Vert^{2}-\rho \Vert z_{k+1}-z_{k+\frac{1}{2}}\Vert^{2}, \notag
\end{align}
where $\rho=\frac{1-\alpha^{2}\kappa^{2}_{m}}{2\alpha}>0$ if $\alpha<\frac{1}{\kappa_{m}}$. Similar to the derivation of \eqref{p7a}, we obtain that $(z_{k+\frac{1}{2}}-z^{*})^{T}F(z_{k+\frac{1}{2}})\ge 0$ for $\forall k\ge 0$. It then follows that
\begin{align*}
\Vert  z_{k+1}-z^{*}\Vert^{2}\le \Vert z_{k}-z^{*}\Vert^{2}-2\alpha\rho \Vert z_{k+1}-z_{k+\frac{1}{2}}\Vert^{2}.
\end{align*}
Similar to the analysis of Theorem 1, we have that $\lim_{k\to \infty} \Vert z_{k}-z^{*}\Vert^{2}=0$. Thus, we conclude that the sequence $\{z_{k}\}$ converges to a saddle point $z^{*}$ of the problem \eqref{1}.

We further analyze the convergence rate of algorithm \eqref{g3}. Let $\hat{x}_{T}=\frac{1}{T}\sum^{T-1}_{k=0}x_{k+\frac{1}{2}}$ and $\hat{y}_{T}=\frac{1}{T}\sum^{T-1}_{k=0}y_{k+\frac{1}{2}}$. From \eqref{lm1} in Lemma 3 and \eqref{p13}, one has that $f(\hat{x}_{T},\lambda^{*})-f(x^{*},\hat{y}_{T}) \notag \le \frac{1}{T}\sum^{T-1}_{k=0}(z_{k+\frac{1}{2}}-z^{*})^{T}F(z_{k+\frac{1}{2}}) \notag \le \frac{1}{2\alpha T}\Vert z_{0}-z^{*}\Vert^{2}$. Similar to the proofs of Theorem 1, we obtain that $\vert f(\hat{x}_{T},\hat{y}_{T})-f(x^{*},y^{*}) \vert \le \frac{1}{2\alpha T}\Vert z_{0}-z^{*}\Vert^{2}$.

%By choosing the parameters , one has that $\frac{1}{2\alpha} \Vert z_{k+1}-z_{k}\Vert^{2}>l_{m}\Vert z_{k+1}-z_{k}\Vert^{2}$.

%\end{IEEEbiography}
%% insert where needed to balance the two columns on the last page with
%%\newpage


%\vfill


\end{document}